\documentclass[11pt,dvips,twoside,letterpaper]{article}
\usepackage[usenames]{color}

\usepackage{pslatex}
\usepackage{fancyhdr}
\usepackage{graphicx}
\usepackage{geometry}

\def\figurename{Figure} 
\makeatletter
\renewcommand{\fnum@figure}[1]{\figurename~\thefigure.}
\makeatother

\def\tablename{Table} 
\makeatletter
\renewcommand{\fnum@table}[1]{\tablename~\thetable.}
\makeatother

\usepackage{amsmath}
\usepackage{amssymb}
\usepackage{amsfonts}
\usepackage{amsthm,amscd}

\newtheorem{theorem}{Theorem}[section]
\newtheorem{lemma}[theorem]{Lemma}
\newtheorem{corollary}[theorem]{Corollary}
\newtheorem{proposition}[theorem]{Proposition}
\theoremstyle{definition}
\newtheorem{definition}[theorem]{Definition}

\newtheorem*{notation}{Notation}

\theoremstyle{remark}
\newtheorem{remark}[theorem]{Remark}

\numberwithin{equation}{section}

\begin{document}
\vskip 0.4in
\title{\bfseries\scshape{Homogenization of Steklov Spectral Problems with Indefinite Density Function in Perforated
Domains}}

\author{\bfseries\scshape Hermann Douanla\thanks{E-mail address: \tt{douanla@chalmers.se}}\\
Department of Mathematical Sciences \\Chalmers University of
Technology\\ Gothenburg, SE-41296, Sweden}

\date{}
\maketitle \thispagestyle{empty}

\begin{abstract} \noindent
The asymptotic behavior of second order self-adjoint elliptic Steklov eigenvalue problems with
periodic rapidly oscillating coefficients and with indefinite (sign-changing) density function is
investigated in periodically perforated domains. We prove that the spectrum of this problem is discrete
 and consists of two sequences, one tending to $-\infty$ and another to $+\infty$. The limiting
 behavior of positive and negative eigencouples  depends crucially on whether the average of the weight
  over the surface of the reference hole is positive,
negative or equal to zero. By means of the two-scale convergence
method, we investigate all three cases.
\end{abstract}


\vspace{.08in}


\section{Introduction}\label{s1}
In 1902, with a motivation coming from Physics, Steklov\cite{Steklov} introduced the following problem
\begin{equation} \label{eq0}
\left\{\begin{aligned} \Delta u&=0\ \ \ \ \text{ in } \Omega\\
\frac{\partial u}{\partial n}&=\rho\lambda u\ \ \ \text{ on }\partial\Omega,
\end{aligned}\right.
\end{equation}
where $\lambda$ is a scalar and $\rho$ is a density function. The
function $u$ represents the steady state temperature on $\Omega$
such that the flux on the boundary $\partial\Omega$ is proportional
to the temperature. In two dimensions, assuming $\rho=1$, problem
(\ref{eq0}) can also be interpreted as a membrane with whole mass
concentrated on the boundary. This problem has been later referred
to as Steklov eigenvalue problem (Steklov is often transliterated as
"Stekloff"). Moreover, eigenvalue problems also arise from many
nonlinear problems after linearization (see e.g., the work of Hess
and Kato\cite{Hess, HessKato} and that of de
Figueiredo\cite{Figueiredo}). This paper deals with the limiting
behavior of a sequence of second order elliptic Steklov eigenvalue
problems with indefinite(sign-changing) density function in
perforated domains.

Let $\Omega$ be a bounded domain in $\mathbb{R}^N_x$(the numerical space of variables $x=(x_1,..., x_N)$),
 with  $\mathcal{C}^1$ boundary $\partial\Omega$ and with integer $N\geq 2$. We define the perforated domain $\Omega^\varepsilon$ as follows.
Let $T\subset Y=(0,1)^N$ be a compact subset of $Y$ with $\mathcal{C}^1$ boundary $\partial T\ (\equiv S)$
and nonempty interior. For $\varepsilon>0$, we define

\[
t^\varepsilon =\{k\in \mathbb{Z}^N: \varepsilon(k+T)\subset\Omega\}
\]
\[
T^\varepsilon=\bigcup_{k\in t^\varepsilon}\varepsilon(k+T)
\]
and
\[
\Omega^\varepsilon=\Omega\setminus T^\varepsilon.
\]
In this setup, $T$ is the reference hole whereas $\varepsilon(k+T)$
is a hole of size $\varepsilon$ and $T^\varepsilon$ is the
collection of the holes of the perforated domain
$\Omega^\varepsilon$. The family $T^\varepsilon$  is made up with a
finite number of holes since $\Omega$ is bounded. In the sequel,
$Y^*$ stands for $Y\setminus T$ and  $n=(n_i)_{i=1}^N$ denotes the
outer unit normal vector to $S$ with respect to $Y^*$.

\pagestyle{fancy} \fancyhead{} \fancyhead[EC]{Hermann Douanla}
\fancyhead[EL,OR]{\thepage} \fancyhead[OC]{Steklov Eigenvalue
Problems with Sing-changing Density Function} \fancyfoot{}
\renewcommand\headrulewidth{0.5pt}

We are interested in the spectral asymptotics (as $\varepsilon\to
0$) of the following Steklov eigenvalue problem
\begin{equation} \label{eq1.1}
\left\{\begin{aligned} -\sum_{i,j=1}^N\frac{\partial}{\partial
x_i}\left(a_{ij}(\frac{x}{\varepsilon})\frac{\partial
u_\varepsilon}{\partial x_j}\right)&=0\text{ in } \Omega^\varepsilon\\
\sum_{i,j=1}^N a_{ij}(\frac{x}{\varepsilon})\frac{\partial
u_\varepsilon}{\partial
x_j}n_i(\frac{x}{\varepsilon})&=\rho(\frac{x}{\varepsilon})\lambda_\varepsilon
u_\varepsilon \text{ on }\partial
T^\varepsilon\\
u_\varepsilon&=0 \text{ on } \partial \Omega,
\end{aligned}\right.
\end{equation}
where $a_{ij}\in L^\infty(\mathbb{R}^N_y)$ ($ 1\leq i,j\leq N$), with the symmetry condition
$a_{ji}=a_{ij}$, the $Y$-periodicity hypothesis: for every
$k\in\mathbb{Z}^N$ one has $a_{ij}(y+k)=a_{ij}(y)$ almost everywhere
in $y\in\mathbb{R}^N_y $, and finally the (uniform) ellipticity condition:
there exists $\alpha>0$ such that
\begin{equation}\label{eq1.2}
 \sum_{i,j=1}^{N}a_{ij}(y)\xi_j\xi_i\geq\alpha|\xi|^2
\end{equation} for all $\xi\in\mathbb{R}^N$ and for almost all $y\in\mathbb{R}^N_y$, where $|\xi|^2=|\xi_1|^2+\cdots +|\xi_N|^2$.
The density function $\rho\in \mathcal{C}_{per}(Y)$ changes sign on $S$, that is, both the set
$\{y\in S, \rho(y)<0\}$ and $\{y\in S, \rho(y)>0\}$ are of
positive $N-1$ dimensional Hausdorf measure (the so-called surface measure). This hypothesis makes the problem under
consideration nonstandard. We will see (Corollary \ref{c11}) that under the preceding
hypotheses, for each $\varepsilon>0$ the spectrum of (\ref{eq1.1})
is discrete and consists of two infinite sequences
$$
0<\lambda_\varepsilon^{1,+}\leq \lambda_\varepsilon^{2,+}\leq \cdots
\leq \lambda_\varepsilon^{n,+}\leq \dots,\quad \lim_{n\to
+\infty}\lambda_\varepsilon^{n,+}=+\infty
$$
and
$$
0>\lambda_\varepsilon^{1,-}\geq \lambda_\varepsilon^{2,-}\geq \cdots
\geq \lambda_\varepsilon^{n,-}\geq \dots,\quad \lim_{n\to
+\infty}\lambda_\varepsilon^{n,-}=-\infty.
$$
The asymptotic behavior of the eigencouples depends crucially on whether the average of the density $\rho$ over $S$,
$M_{S}(\rho)=\int_{S}\rho(y)\,d\sigma(y)$, is positive, negative or equal to zero. All three cases are carefully investigated in this paper.

The homogenization of spectral problems has been widely explored. In
a fixed domain, homogenization of spectral problems with point-wise
positive density function goes back to Kesavan \cite{Kesavan1,
Kesavan2}. Spectral asymptotics in perforated domains was studied by
Vanninathan\cite{Vanni} and later in many other papers, including
\cite{douanla1, douanla2, Kaizu,  oleinik, Pastukhova, Roppongi2}
and the references therein. Homogenization of elliptic operators
with sing-changing density function in a fixed domain with Dirichlet
boundary conditions has been investigated by Nazarov et al.
\cite{Nazarov1, Nazarov2, Nazarov3} via a combination of formal
asymptotic expansion with Tartar's energy method. In porous media,
spectral asymptotics of elliptic operator with sign changing density
function is studied in \cite{douanla3} with the two scale
convergence method.

The asymptotics of Steklov eigenvalue problems in periodically
perforated domains was studied in \cite{Vanni} for the laplace
operator and constant density ($\rho=1$) using asymptotic expansion
and Tartar's test function method. The same problem for a second
order periodic elliptic operator has been studied in
\cite{Pastukhova} (with $\mathcal{C}^\infty$ coefficients) and in
\cite{douanla1} (with $L^\infty$ coefficient) but still with
constant density ($\rho=1$). All the just-cited works deal only with
one sequence of positive eigenvalues.

In this paper we take it to the general tricky step. We investigate
in periodically perforated domains the asymptotic behavior of
Steklov eigenvalue problems for periodic elliptic linear
differential operators of order two in divergence form with
$L^\infty$ coefficients and a sing-changing density function. We
obtain accurate and concise homogenization results in all three
cases: $ M_{S}(\rho)>0 $ (Theorem \ref{t3.1} and Theorem
\ref{t3.2}),\ $M_{S}(\rho)=0$ (Theorem \ref{t3.3}),\ $M_{S}(\rho)<0$
(Theorem \ref{t3.1}
 and Theorem \ref{t3.2}), by using the two-scale
convergence method\cite{AB, GNWL, G89, Zhikov} introduced by
Nguetseng\cite{G89} and further developed by Allaire\cite{AB}. In short;
\begin{description}
  \item[i)] If $M_{S}(\rho)>0$, then the positive eigencouples behave like in the case of point-wise positive density function, i.e.,
for $k\geq 1$, $\lambda^{k,+}_\varepsilon$ is of order $\varepsilon$  and $\frac{1}{\varepsilon}\lambda^{k,+}_\varepsilon$ converges as $\varepsilon\to 0$ to the $k^{th}$ eigenvalue of the limit
 Dirichlet spectral problem, corresponding extended eigenfunctions converge along subsequences.

 As regards the "negative" eigencouples, $\lambda^{k,-}_\varepsilon$ converges to $-\infty$ at the rate $\frac{1}{\varepsilon}$ and
 the corresponding eigenfunctions oscillate rapidly. We use a factorization technique (\cite{Kozlov, Vanni}) to prove that
 $$
 \lambda^{k,-}_\varepsilon=\frac{1}{\varepsilon}\lambda_1^-+\xi^{k,-}_\varepsilon +o(1),\qquad k=1,2\cdots
 $$
 where
   $(\lambda_1^-, \theta_{1}^-)$
    is the first negative eigencouple to the following local Steklov spectral problem
    \begin{equation}
\left\{\begin{aligned}-div(a(y)D_y\theta)&=0 \quad\text{ in } Y^*\\
a(y)D_y\theta\cdot n&=\lambda\rho(y)\theta \quad\text{ on }\ \ S\\
\theta \quad &Y-periodic,
\end{aligned}\right.
\end{equation}
   and $\{\xi^{k,\pm}_\varepsilon\}_{k=1}^\infty$ are eigenvalues of a Steklov eigenvalue problem similar  to (\ref{eq1.1}). We then prove that $\{\frac{\lambda^{k,-}_\varepsilon}{\varepsilon}-\frac{\lambda_1^-}{\varepsilon^2}\}$ converges to the $k^{th}$ eigenvalue of  a limit Dirichlet spectral problem
     which is different from that obtained for positive eigenvalues. As regards eigenfunctions,
    extensions of  $\{\frac{u^{k,-}_\varepsilon}{(\theta_{1}^-)^\varepsilon}\}_{\varepsilon\in E}$ - where  $(\theta_{1}^-)^\varepsilon(x)=\theta_{1}^-(\frac{x}{\varepsilon})$ -
      converge along subsequences to the $k^{th}$ eigenfunctions of the limit problem.
 \item[ii)] If $M_{S}(\rho)=0$, then the limit spectral problem generates a quadratic operator pencil and
 $\lambda^{k,\pm}_\varepsilon$ converges to the $(k,\pm)^{th}$ eigenvalue of the limit operator, extended eigenfunctions converge along
 subsequences as well. This case requires a new convergence result as regards the two-scale convergence theory, Lemma \ref{l2.1}.
\item[iii)] The case when $M_{S}(\rho)<0$ is equivalent to that when $M_{S}(\rho)>0$, just replace $\rho$ with $-\rho$.
\end{description}
Unless otherwise specified, vector spaces throughout are considered
over $\mathbb{R}$, and scalar functions are assumed to take real
values. We will make use of the following notations. Let
$F(\mathbb{R}^N)$ be a given function space. We denote by
$F_{per}(Y)$ the space of functions in $F_{loc}(\mathbb{R}^N)$ (when
it makes sense) that are $Y$-periodic, and by
$F_{per}(Y)/\mathbb{R}$ the space of those functions $u\in
F_{per}(Y)$ with $\int_Y u(y)dy=0$. We denote by $H^1_{per}(Y^*)$
the space of functions in $H^1(Y^*)$ assuming same values on the
opposite faces of $Y$ and $H^1_{per}(Y^*)/\mathbb{R}$ stands for the
subset of $H^1_{per}(Y^*)$ made up of functions $u\in
H^1_{per}(Y^*)$ verifying $\int_{Y^*}u(y)dy=0$. Finally, the letter
$E$ denotes throughout a family of  strictly positive real numbers
$(0<\varepsilon<1)$ admitting $0$ as accumulation point. The
numerical space $\mathbb{R}^N$ and its open sets are provided with
the Lebesgue measure denoted by $dx=dx_1...dx_N$. The usual gradient
operator will be denoted by $D$. For the sake of simple notations we
hide trace operators. The rest of the paper is organized as follows.
Section \ref{s2} deals with some preliminary results while
homogenization processes are considered in Section \ref{s3}.

\section{Preliminaries}\label{s2}
We first recall the definition and the main compactness theorems of
the two-scale convergence method. Let $\Omega$ be a smooth open
bounded set in $\Bbb{R}^N_x$ (integer $N\geq 2$) and $Y=(0,1)^N$,
the unit cube.

\begin{definition}\label{d2.1}
A sequence $(u_\varepsilon)_{\varepsilon\in E}\subset L^2(\Omega)$
is said to two-scale converge in $L^2(\Omega)$ to some $u_0\in
L^2(\Omega\times Y)$ if as $E\ni\varepsilon\to 0$,
\begin{equation}\label{eq2.14}
    \int_\Omega u_\varepsilon(x)\phi(x,\frac{x}{\varepsilon})dx\to \iint_{\Omega\times Y}u_0(x,y)\phi(x,y)dxdy
\end{equation}
for all $\phi\in L^2(\Omega;\mathcal{C}_{per}(Y))$.
\end{definition}
\begin{notation}
We express this by writing $u_\varepsilon\xrightarrow{2s}u_0$ in
$L^2(\Omega)$.
\end{notation}
The following compactness theorems \cite{AB, G89, GW2007} are cornerstones of the two-scale convergence
method.
\begin{theorem}\label{t2.1}
Let $(u_\varepsilon)_{\varepsilon\in E}$ be a bounded sequence in
$L^2(\Omega)$. Then a subsequence $E'$ can be extracted from  $E$
such that as $E'\ni\varepsilon\to 0$, the sequence
$(u_\varepsilon)_{\varepsilon\in E'}$ two-scale converges in
$L^2(\Omega)$ to some $u_0\in L^2(\Omega\times Y)$.
\end{theorem}

\begin{theorem}\label{t2.21}
Let $(u_\varepsilon)_{\varepsilon\in E}$ be a bounded sequence in
$H^1(\Omega)$. Then a subsequence $E'$ can be extracted from  $E$
such that as $E'\ni\varepsilon\to 0$
\begin{eqnarray}
  u_\varepsilon &\to& u_0 \ \ \  \text{ in } H^1(\Omega)\text{-weak}\nonumber\\
  u_\varepsilon &\to& u_0  \ \ \ \ \ \ \ \ \text{ in } L^2(\Omega)\nonumber\\
  \frac{\partial u_\varepsilon}{\partial x_j}&\xrightarrow{2s}&\frac{\partial u_0}{\partial x_j}+\frac{\partial u_1}{\partial y_j}
  \ \ \ \ \text{ in } L^2(\Omega) \ \ (1\leq j \leq N)\nonumber
\end{eqnarray}
where $u_0\in H^1(\Omega)$ and $u_1\in L^2(\Omega;H^1_{per}(Y))$.
Moreover, as $E'\ni\varepsilon\to 0$ we have
\begin{equation}\label{eq1.14}
\int_{\Omega}\frac{u_\varepsilon (x)}{\varepsilon}\psi(x,\frac{x}{\varepsilon})dx\to
\iint_{\Omega\times Y}u_1(x,y)\psi(x,y)dx\,dy
\end{equation}
for $\psi\in \mathcal{D}(\Omega)\otimes (L^2_{per}(Y)/\mathbb{R})$.
\end{theorem}

\begin{remark}
In Theorem \ref{t2.21} the function $u_1$ is unique up to an additive function of variable $x$. We need to fix its choice according to our
future needs. To do this, we introduce the following space
\[
H^{1,*}_{per}(Y)=\{u\in H^1_{per}(Y):\int_{Y^*} u(y)dy=0\}.
\]
\noindent This defines a closed subspace of $H^1_{per}(Y)$ as it is the kernel of the bounded linear
 functional $u\mapsto\int_{Y^*}u(y)dy$ defined on $H^1_{per}(Y)$. It is to be noted that for $u\in H^{1,*}_{per}(Y)$,
 its restriction to $Y^*$ (which will still be denoted by $u$ in the sequel) belongs to $H^1_{per}(Y^*)/\mathbb{R}$.
\end{remark}
 We will use the following version of Theorem \ref{t2.21}.

\begin{theorem}\label{t2.2}
Let $(u_\varepsilon)_{\varepsilon\in E}$ be a bounded sequence in
$H^1(\Omega)$. Then a subsequence $E'$ can be extracted from  $E$
such that as $E'\ni\varepsilon\to 0$
\begin{eqnarray}
  u_\varepsilon &\to& u_0 \ \ \  \text{ in } H^1(\Omega)\text{-weak}\label{eq1.111}\\
  u_\varepsilon &\to& u_0  \ \ \ \ \ \ \ \ \text{ in } L^2(\Omega)\label{eq1.121}\\
  \frac{\partial u_\varepsilon}{\partial x_j}&\xrightarrow{2s}&\frac{\partial u_0}{\partial x_j}+\frac{\partial u_1}{\partial y_j}
  \ \ \ \ \text{ in } L^2(\Omega) \ \ (1\leq j \leq N)\label{eq1.131}
\end{eqnarray}
where $u_0\in H^1(\Omega)$ and $u_1\in L^2(\Omega;H^{1,*}_{per}(Y))$.
Moreover, as $E'\ni\varepsilon\to 0$ we have
\begin{equation}\label{eq1.1311}
\int_{\Omega}\frac{u_\varepsilon (x)}{\varepsilon}\psi(x,\frac{x}{\varepsilon})dx\to
\iint_{\Omega\times Y}u_1(x,y)\psi(x,y)dx\,dy
\end{equation}
for $\psi\in \mathcal{D}(\Omega)\otimes (L^2_{per}(Y)/\mathbb{R})$.
\end{theorem}
\begin{proof}
Let $\widetilde{u}_1\in L^2(\Omega; H^1_{per}(Y))$ be such that Theorem \ref{t2.21} holds with $\widetilde{u}_1$ in place of $u_1$. Put
\[
u_1(x,y)=\widetilde{u}_1(x,y)-\frac{1}{|Y^*|}\int_{Y^*}\widetilde{u}_1(x,y)dy\qquad (x,y)\in\Omega\times Y,
\]
where $|Y^*|$ stands for the Lebesgue measure of $Y^*$. Then $u_1\in L^2(\Omega; H^{1,*}_{per}(Y))$ and moreover
 $D_y u_1= D_y\widetilde{u}_1$ so that (\ref{eq1.131}) holds.
\end{proof}

In the sequel, $S^\varepsilon$ stands for $\partial T^\varepsilon$
and the surface measures on $S$ and $S^\varepsilon$ are denoted by
$d\sigma(y)$ ($y\in Y$), $d\sigma_\varepsilon(x)$ ($x\in\Omega,
\varepsilon\in E$), respectively. The space of squared integrable
functions, with respect to the previous measures on $S$ and
$S^\varepsilon$ are denoted by $L^2(S)$ and $L^2(S^\varepsilon)$
respectively. Since the volume of $S^\varepsilon$ grows
proportionally to $\frac{1}{\varepsilon}$ as $\varepsilon\to 0$, we
endow $L^2(S^\varepsilon)$ with the scaled scalar product\cite{ADH,Radu, Radu2}
$$
\left(u,v\right)_{L^2(S^\varepsilon)}=\varepsilon\int_{S^\varepsilon}u(x)v(x)d\sigma_\varepsilon(x)\quad \left(u,v\in L^2(S^\varepsilon)\right).
$$
Definition \ref{d2.1} and  theorem \ref{t2.1} then generalize as
\begin{definition}\label{d2.2}
A sequence $(u_\varepsilon)_{\varepsilon\in E}\subset
L^2(S^\varepsilon)$ is said to two-scale converge to some $u_0\in
L^2(\Omega\times S)$ if as $E\ni\varepsilon\to 0$,
$$
\varepsilon\int_{S^\varepsilon}
u_\varepsilon(x)\phi(x,\frac{x}{\varepsilon})d\sigma_\varepsilon(x)\to
\iint_{\Omega\times S}u_0(x,y)\phi(x,y)dxd\sigma(y)
$$
for all $\phi\in \mathcal{C}(\overline{\Omega};\mathcal{C}_{per}(Y))$.
\end{definition}
\begin{theorem}\label{t2.3}
Let $(u_\varepsilon)_{\varepsilon\in E}$ be a sequence in $L^2(S^\varepsilon)$ such that
$$
\varepsilon\int_{S^\varepsilon}|u_\varepsilon(x)|^2d\sigma_\varepsilon(x)\leq C
$$
where $C$ is a positive constant independent of $\varepsilon$. There
exists a subsequence $E'$ of $E$ such that
$(u_\varepsilon)_{\varepsilon\in E'}$ two-scale converges to some
$u_0\in L^2(\Omega;L^2(S))$ in the sense of definition \ref{d2.2}.
\end{theorem}

In the case when $(u_\varepsilon)_{\varepsilon\in E}$ is the
sequence of traces on $S^\varepsilon$ of functions in $H^1(\Omega)$,
one can link its usual two-scale limit with its surface two-scale
limits. The following proposition whose proof can be found in
\cite{ADH} clarifies this.

\begin{proposition}\label{p2.1}
Let $(u_\varepsilon)_{\varepsilon\in E}\subset H^1(\Omega)$ be such that
\begin{equation*}\label{}
    \|u_\varepsilon\|_{L^2(\Omega)}+\varepsilon\|Du_\varepsilon\|_{L^2(\Omega)^N}\leq C,
\end{equation*}
where $C$ is a positive constant independent of $\varepsilon$ and
$D$ denotes the usual gradient. The sequence of traces of
$(u_\varepsilon)_{\varepsilon\in E}$ on $S^\varepsilon$ satisfies
$$
\varepsilon\int_{S^\varepsilon}|u_\varepsilon(x)|^2d\sigma_\varepsilon(x)\leq
C\quad (\varepsilon\in E)
$$
and up to a subsequence $E'$ of $E$, it two-scale converges in the
sense of Definition \ref{d2.2} to some $u_0\in L^2(\Omega;L^2(S))$
 which is nothing but the trace on $S$ of the usual two-scale limit, a function in $L^2(\Omega;H^1_{per}(Y))$. More precisely,
 as $E'\ni\varepsilon\to 0$

 \begin{eqnarray*}
    \varepsilon\int_{S^\varepsilon}u_\varepsilon(x)\phi(x,\frac{x}{\varepsilon})d\sigma_\varepsilon(x)&\to&\iint_{\Omega\times S}
    u_0(x,y)\phi(x,y)dxd\sigma(y),\\
    \int_{\Omega}u_\varepsilon(x)\phi(x,\frac{x}{\varepsilon})dxdy&\to&\iint_{\Omega\times Y}u_0(x,y)\phi(x,y)dxdy,
 \end{eqnarray*}
 for all $\phi\in\mathcal{C}(\overline{\Omega};\mathcal{C}_{per}(Y))$.
\end{proposition}

In our homogenization process,  more precisely in the case when $M_S(\rho)=0$, we will need a generalization of
(\ref{eq1.14}) to periodic surfaces. Notice that (\ref{eq1.14}) was proved for the first time in a deterministic
setting by Nguetseng and Woukeng in \cite{GW2007} but to the best of our knowledge its generalization to periodic
surfaces is not yet available in the literature. We state and prove it below.
\begin{lemma}\label{l2.1}
Let $(u_\varepsilon)_{\varepsilon\in E}\subset H^1(\Omega)$ be such that as $E\ni \varepsilon\to 0$
\begin{eqnarray}
  u_\varepsilon &\xrightarrow{2s}& u_0 \ \ \  \text{ in } L^2(\Omega)\label{eq1.11}\\
   \frac{\partial u_\varepsilon}{\partial x_j}&\xrightarrow{2s}&\frac{\partial u_0}{\partial x_j}+\frac{\partial u_1}{\partial y_j}
  \ \ \ \ \text{ in } L^2(\Omega) \ \ (1\leq j \leq N)\label{eq1.13}
\end{eqnarray}
for some $u_0\in H^1(\Omega)$ and $u_1\in L^2(\Omega; H^1_{per}(Y))$. Then
\begin{equation}\label{}
    \lim_{\varepsilon\to 0}\int_{S^\varepsilon}u_\varepsilon(x)\varphi(x)\theta(\frac{x}{\varepsilon})d\sigma_\varepsilon(x)=\iint_{\Omega\times S}u_1(x,y)\varphi(x)\theta(y)dxd\sigma(y)
\end{equation}
for all $\varphi\in\mathcal{D}(\Omega)$ and $\theta\in\mathcal{C}_{per}(Y)$ with $\int_S \theta(y)d\sigma(y)=0$.
\end{lemma}
\begin{proof}
We first transform the above surface integral into a volume integral by adapting the trick in \cite[Section 3]{donato}.
 By the mean value zero condition over $S$ for $\theta$ we conclude that there exists
    a unique solution $\vartheta\in H^1_{per}(Y^*)/\mathbb{R}$ to
\begin{equation}\label{eq3.998}
\left\{\begin{aligned} &-\Delta_y \vartheta=0\ \ \text{in } Y^*\\
&D_y\vartheta(y)\cdot n(y)=\theta(y)\ \ \text{on } S,
\end{aligned}\right.
\end{equation}
where $n=(n_i)_{i=1}^N$ stands for the outward unit normal to $S$ with respect to $Y^*$. Put $\phi=D_y \vartheta$. We get
\begin{eqnarray}
  &&\int_{\Omega^\varepsilon}D_x u_\varepsilon(x)\varphi(x)\cdot D_y \vartheta(\frac{x}{\varepsilon})dx =
  \int_{S^\varepsilon}u_\varepsilon(x)\varphi(x)D_y\vartheta(\frac{x}{\varepsilon})\cdot n(\frac{x}{\varepsilon})d\sigma_\varepsilon(x)\nonumber \\
  & &\qquad\qquad\quad\qquad- \int_{\Omega^\varepsilon}u_\varepsilon(x)D_x\varphi(x)\cdot D_y\vartheta(\frac{x}{\varepsilon})dx
  -\frac{1}{\varepsilon}\int_{\Omega^\varepsilon}u_\varepsilon(x)\varphi(x)\Delta_y \vartheta(\frac{x}{\varepsilon})dx \qquad\qquad\label{eq3.999}\\
   &&\qquad\qquad\qquad=\int_{S^\varepsilon}u_\varepsilon(x)\varphi(x)\theta(\frac{x}{\varepsilon})d\sigma_\varepsilon(x)-
   \int_{\Omega^\varepsilon}u_\varepsilon(x)D_x\varphi(x)\cdot \phi(\frac{x}{\varepsilon})dx.\nonumber
\end{eqnarray}
Next, sending $\varepsilon$ to $0$ yields
\begin{eqnarray*}
    \lim_{\varepsilon\to 0}\int_{S^\varepsilon}u_\varepsilon(x)\varphi(x)\theta(\frac{x}{\varepsilon})\,d\sigma_\varepsilon(x)&=& \iint_{\Omega\times Y^*}
    [D_x u_0(x)+D_y u_1(x,y)]\varphi(x)\cdot \phi(y)\,dxdy \\
    & &+\iint_{\Omega\times Y^*}u_0(x)D_x\varphi(x)\cdot \phi(y)\,dxdy \\
   &=& \iint_{\Omega\times Y^*}D_y u_1(x,y)\varphi(x)\cdot\phi(y)\,dxdy.
\end{eqnarray*}
We finally have
\begin{eqnarray*}
   \iint_{\Omega\times Y^*}D_y u_1(x,y)\varphi(x)\cdot\phi(y)\,dxdy  &=& -\iint_{\Omega\times Y^*}u_1(x,y)\varphi(x)\Delta_y\vartheta(y)\,dxdy \\
   &+&\iint_{\Omega\times S} u_1(x,y)\varphi(x)\phi(y)\cdot n(y)\,dxd\sigma(y) \\
  &=&  \iint_{\Omega\times S} u_1(x,y)\varphi(x)\theta(y)\,dxd\sigma(y).
\end{eqnarray*}
The proof is completed.
\end{proof}

We now gather some preliminary results. We introduce the
characteristic function $\chi_{G}$ of
$$
G=\mathbb{R}^N_y\setminus \Theta
$$
with
$$
\Theta=\bigcup_{k\in \mathbb{Z}^N}(k+T).
$$
It is clear that $G$ is an open subset of $\mathbb{R}^N_y$. Next,
let $\varepsilon\in E$ be arbitrarily fixed and define
\[
V_\varepsilon =\{u\in H^1(\Omega^\varepsilon) : u=0 \text{ on }
\partial\Omega\}.
\]
We equip $V_\varepsilon$ with the $H^1(\Omega^\varepsilon)$-norm
which makes it a Hilbert space. We recall the following classical
extension result \cite{CSJP}.
\begin{proposition}\label{p3.1}
For each $\varepsilon\in E$ there exists an operator $P_\varepsilon$
of  $V_\varepsilon$ into $H^1_0(\Omega)$ with the following
properties:
\begin{itemize}
    \item $P_\varepsilon$ sends continuously and linearly $V_\varepsilon$ into $H^1_0(\Omega)$.
    \item $(P_\varepsilon v)|_{\Omega^\varepsilon}=v$ for all $v\in V_\varepsilon$.
    \item $\|D(P_\varepsilon v)\|_{L^2(\Omega)^N}\leq c \|D v\|_{L^2(\Omega^\varepsilon)^N}$ for all $v\in V_\varepsilon$,
     where $c$ is a constant independent of $\varepsilon$.
\end{itemize}
\end{proposition}
In the sequel, we will explicitly write  the just-defined extension operator everywhere needed but we will abuse notations
on the local extension operator (see \cite{CSJP} for its definition): the extension to $Y$ of $u\in H^1_{per}(Y^*)/\mathbb{R}$
will still be denoted by $u$ (this extension is an element of $H^{1,*}_{per}(Y)$).

Now, let $Q^\varepsilon=\Omega\setminus (\varepsilon \Theta)$. This
defines an open set in $\mathbb{R}^N$ and
$\Omega^\varepsilon\setminus Q^\varepsilon$ is the intersection of
$\Omega$ with the collection of the holes crossing the boundary
$\partial\Omega$. The following result implies that the holes
crossing the boundary $\partial\Omega$ are of no effects as regards
the homogenization process.

 \begin{lemma}\cite{Gper}\label{l3.1}
 Let $K\subset\Omega$ be a compact set independent of $\varepsilon$. There is some $\varepsilon_0>0$ such
  that $\Omega^\varepsilon\setminus Q^\varepsilon\subset\Omega\setminus K$ for any $0<\varepsilon\leq\varepsilon_0$.
 \end{lemma}

We introduce the space
\[
\mathbb{F}^1_0=H^1_0(\Omega)\times L^2\left(\Omega; H^{1,*}_{per}(Y)\right)
\]
and endow it with the following norm
\[
\|\textbf{v}\|_{\mathbb{F}^1_0}= \left\|D_x v_0+D_y
v_1\right\|_{L^2(\Omega\times Y)}\ \ \ \
(\textbf{v}=(v_0,v_1)\in\mathbb{F}^1_0),
\]
which makes it a Hilbert space admitting $F_{0}^{\infty
}=\mathcal{D}(\Omega )\times [\mathcal{D}(\Omega )\otimes
\mathcal{C}_{per}^{\infty,*}(Y)]$ (where
$\mathcal{C}_{per}^{\infty,*}(Y)=\{u\in\mathcal{C}_{per}^{\infty}(Y):\int_{Y^*}u(y)dy=0\}$)
as a dense subspace. For
  $(\textbf{u},\textbf{v})\in\mathbb{F}^1_0\times\mathbb{F}^1_0$, let
\begin{equation*}
    a_\Omega(\textbf{u},\textbf{v})=\sum_{i,j=1}^N\iint_{\Omega \times Y^*}a_{ij}(y)\left(\frac{\partial u_0}{\partial x_j}+
    \frac{\partial u_1}{\partial y_j}\right)\left(\frac{\partial v_0}{\partial x_i}+
    \frac{\partial v_1}{\partial y_i}\right)\,dxdy.
 \end{equation*}
This define a symmetric, continuous bilinear form on
$\mathbb{F}^1_0\times \mathbb{F}^1_0$. We will need the following
results whose proof can be found in \cite{douanla1}.
\begin{lemma}\label{l3.2}
Fix $\Phi=(\psi _{0},\psi _{1})\in F_0^\infty$ and define
$\Phi_\varepsilon:\Omega\to \mathbb{R}$ ($\varepsilon>0$) by
\begin{equation*}
    \Phi_\varepsilon(x)=\psi_0(x) +
    \varepsilon\psi_1(x,\frac{x}{\varepsilon})\quad (x\in \Omega).
\end{equation*}
If $(u_\varepsilon)_{\varepsilon\in E}\subset H^1_0(\Omega)$ is such
that
\begin{equation*}
    \frac{\partial u_\varepsilon}{\partial x_i}\xrightarrow{2s} \frac{\partial u_0}{\partial x_i}+\frac{\partial u_1}
    {\partial y_i}\ \ \text{ in }\ \ L^2(\Omega) \ (1\leq i\leq N)
\end{equation*}
as $E\ni \varepsilon\to 0$ for some $\textbf{u}=(u_0, u_1)\in
\mathbb{F}^1_0$, then
 \begin{equation*}
    a^\varepsilon(u_\varepsilon,\Phi_\varepsilon)\to a_\Omega(\textbf{u},\Phi)
 \end{equation*}
as  $E\ni \varepsilon\to 0$, where
\[
 a^\varepsilon(u_\varepsilon,\Phi_\varepsilon)=\sum_{i,j=1}^N\int_{\Omega^\varepsilon} a_{ij}(\frac{x}{\varepsilon})\frac{\partial u_\varepsilon}{\partial x_j}
\frac{ \partial \Phi_\varepsilon}{\partial x_i}dx.
\]
\end{lemma}

We now construct and point out the main properties of the so-called
homogenized coefficients. Put
\begin{equation*}
    a(u,v)=\sum_{i,j=1}^N\int_{Y^*}a_{ij}(y)\frac{\partial u}{\partial y_j}\frac{\partial v}{\partial
    y_i}dy,
\end{equation*}
\begin{equation*}
   \qquad\qquad l_j(v)=\sum_{k=1}^N\int_{Y^*}a_{kj}(y)\frac{\partial v}{\partial
    y_k}dy,\qquad (1\leq j\leq N)
\end{equation*}
and
\begin{equation*}
    l_0(v)=\int_{S}\rho(y)v(y)d\sigma(y),
\end{equation*}
 for $u,v\in H^1_{per}(Y^*)/\mathbb{R}$. Equipped with the norm
\begin{equation}\label{eq3.311}
   \|u\|_{H^1_{per}(Y^*)/\mathbb{R}}=\|D_y u\|_{L^2(Y^*)^N}\ \ \ (u\in H^1_{per}(Y^*)/\mathbb{R}),
\end{equation}
$H^1_{per}(Y^*)/\mathbb{R}$ is a Hilbert space.

\begin{proposition}\label{p3.2}
Let $1\leq j\leq N$.  The local variational problems
\begin{equation}\label{eq3.9}
  u\in H^1_{per}(Y^*)/\mathbb{R} \text{ and } a(u,v)=l_j(v)\  \text{ for all }\  v\in H^1_{per}(Y^*)/\mathbb{R}
\end{equation}
and
\begin{equation}\label{eq3.91}
  u\in H^1_{per}(Y^*)/\mathbb{R} \text{ and } a(u,v)=l_0(v)\  \text{ for all }\
v\in H^1_{per}(Y^*)/\mathbb{R}
\end{equation}
admit each a unique solution, assuming for (\ref{eq3.91}) that $M_S(\rho)=0$.
\end{proposition}
Let $1\leq i,j\leq N$. The homogenized coefficients read
\begin{equation}\label{eq3.10}
    q_{ij}=\int_{Y^*}a_{ij}(y)dy-\sum_{l=1}^N\int_{Y^*}a_{il}(y)\frac{\partial\chi^j}{\partial
    y_l}(y)dy
\end{equation}
where $\chi^j\ \ (1\leq j\leq N)$ is the solution to (\ref{eq3.9}).
We recall that
 $q_{ji}=q_{ij}\ \ (1\leq i,j\leq N)$ and there exists a constant $\alpha_{0}>0$ such that
 \[
   \sum_{i,j=1}^{N}q_{ij}\xi_{j}\xi_{i}\geq\alpha_{0}|\xi|^{2}
\]
for all $\xi\in \mathbb{R}^{N}$ (see e.g., \cite{BLP}).

We now visit the existence result for (\ref{eq1.1}). The weak formulation of (\ref{eq1.1}) reads: Find
 $(\lambda_\varepsilon, u_\varepsilon)\in\mathbb{C}\times V_\varepsilon$, ($u_\varepsilon\neq 0$) such that
\begin{equation}\label{eq3.101}
    a^\varepsilon(u_\varepsilon,v)=\lambda_\varepsilon(\rho^\varepsilon u_\varepsilon, v)_{S^\varepsilon}, \quad v\in V_\varepsilon,
\end{equation}
where
$$
(\rho^\varepsilon u_\varepsilon, v)_{S^\varepsilon}=\int_{S^\varepsilon}\rho^\varepsilon u_\varepsilon v d\sigma_\varepsilon(x).
$$
Since $\rho^\varepsilon$ changes sign, the classical results on the spectrum of semi-bounded self-adjoint operators with compact
resolvent do not apply. To handle this, we follow the ideas in \cite{Nazarov3}. The bilinear form
$(\rho^\varepsilon u, v)_{S^\varepsilon}$ defines a bounded linear operator $K^\varepsilon:V_\varepsilon\to V_\varepsilon$ such that
$$
(\rho^\varepsilon u, v)_{S^\varepsilon}=a^\varepsilon(K^\varepsilon u,v) \quad (u,v\in V_\varepsilon).
$$
The operator $K^\varepsilon$ is symmetric and its domains
$D(K^\varepsilon)$ coincides with the whole $V_\varepsilon$, thus it
is self-adjoint.  Recall that the gradient norm is equivalent to the
$H^1(\Omega^\varepsilon)$-norm on $V_\varepsilon$. Looking at
$K^\varepsilon u$ as the solution to the boundary value problem
\begin{equation}\label{}\left\{\begin{aligned}
 -div(a(\frac{x}{\varepsilon})D_x(K^\varepsilon u))&=0\quad \text{in }\Omega^\varepsilon\\
 a(\frac{x}{\varepsilon})D_xK^\varepsilon u\cdot n(\frac{x}{\varepsilon})&=\rho^\varepsilon u \quad \text{ on  } S^\varepsilon\\
 K^\varepsilon u(x)&=0 \quad \text{ on  } \partial \Omega,
\end{aligned}\right.
\end{equation}
we get a constant $C_\varepsilon>0$ such that $\|K^\varepsilon
u\|_{V^\varepsilon}\leq C_\varepsilon\|u\|_{L^2(S^\varepsilon)}$.
But the trace operator $V_\varepsilon\to L^2(S^\varepsilon)$ is
compact. The compactness of $K^\varepsilon$ follows thereby. We can
rewrite (\ref{eq3.101}) as follows
$$
K^\varepsilon u_\varepsilon=\mu_\varepsilon u_\varepsilon,\quad \mu_\varepsilon=\frac{1}{\lambda_\varepsilon}.
$$
We recall that (see e.g., \cite{Bs3}) in the case $\rho\geq 0$ on
$S$, the operator $K^\varepsilon$ is positive and its spectrum
$\sigma(K^\varepsilon)$ lies in $[0, \|K^\varepsilon\|]$ and
$\mu_\varepsilon=0$ belongs to the essential spectrum
$\sigma_e(K^\varepsilon)$. Let $L$ be a self-adjoint operator and
let  $\sigma_p^\infty(L)$ and $\sigma_c(L)$ be its set of
eigenvalues of infinite multiplicity and its continuous spectrum,
respectively. We have
$\sigma_e(L)=\sigma_p^\infty(L)\cup\sigma_c(L)$ by definition. The
spectrum of $K^\varepsilon$ is described by the following
proposition whose proof is  similar to that of \cite[Lemma
1]{Nazarov3}.
\begin{lemma}\label{l2.2}
Let $\rho\in\mathcal{C}_{per}(Y)$ be such that the sets $\{y\in S: \rho(y)< 0\}$ and $\{y\in S: \rho(y)> 0\}$ are
both of positive surface measure. Then for any $\varepsilon>0$, we have $\sigma(K^\varepsilon)\subset [-\|K^\varepsilon\|, \|K^\varepsilon\|]$
and $\mu=0$ is the only element of the essential spectrum $\sigma_e(K^\varepsilon)$. Moreover, the discrete spectrum
of $K^\varepsilon$ consists of two infinite sequences
\begin{eqnarray*}
  \mu_\varepsilon^{1,+}\geq \mu_\varepsilon^{2,+}\geq \cdots \geq\mu_\varepsilon^{k,+}\geq\cdots\to 0^+,&& \\
 \mu_\varepsilon^{1,-}\leq \mu_\varepsilon^{2,-}\leq \cdots \leq\mu_\varepsilon^{k,-}\leq\cdots\to 0^-.&& \\
 \end{eqnarray*}
\end{lemma}
\begin{corollary}\label{c11}
The hypotheses are those of Lemma \ref{l2.2}. Problem (\ref{eq1.1}) has a discrete set of eigenvalues consisting of two sequences
\begin{eqnarray*}
  0<\lambda_\varepsilon^{1,+}\leq\lambda_\varepsilon^{2,+}\leq \cdots\leq \lambda_\varepsilon^{k,+}\leq\cdots \to +\infty, && \\
   0>\lambda_\varepsilon^{1,+}\geq\lambda_\varepsilon^{2,-}\geq \cdots\geq \lambda_\varepsilon^{k,-}\geq\cdots \to -\infty. &&
\end{eqnarray*}
\end{corollary}
We may now address the homogenization problem for (\ref{eq1.1}).

\section{Homogenization results}\label{s3}
In this section we state and prove homogenization results for both cases $M_{S}(\rho)>0$ and  $M_{S}(\rho)=0$.
The homogenization results in the case when $M_{S}(\rho)<0$ can be deduced from the case $M_{S}(\rho)>0$ by replacing $\rho$ with $-\rho$.
We start with the less technical case.
\subsection{The case $M_{S}(\rho)>0$}
We start with the homogenization result for the positive part of the
spectrum $(\lambda_\varepsilon^{k,+}, u_\varepsilon^{k,+})_{\varepsilon\in E}$.
\subsubsection{Positive part of the spectrum}
We assume (this is not a restriction) that the corresponding eigenfunctions
are orthonormalized as follows
\begin{equation}\label{eq3.131}
    \varepsilon\int_{S^\varepsilon}\rho(\frac{x}{\varepsilon})u_\varepsilon^{k,+}u_\varepsilon^{l,+}d\sigma_\varepsilon(x)=\delta_{k,l}\quad
    k,l=1,2,\cdots
\end{equation}
and the homogenization results states as

\begin{theorem}\label{t3.1}
For each $k\geq 1$ and each $\varepsilon\in E$, let
$(\lambda^{k,+}_\varepsilon,u^{k,+}_\varepsilon)$ be the $k^{th}$
positive eigencouple to (\ref{eq1.1}) with $M_{S}(\rho)>0$ and
(\ref{eq3.131}). Then, there exists a subsequence $E'$ of $E$ such
that
\begin{eqnarray}
 \frac{1}{\varepsilon}\lambda^{k,+}_\varepsilon &\to&  \lambda^{k}_0\quad\text{in }\ \mathbb{R}\ \text{ as }E\ni\varepsilon\to 0\label{eq3.14}\\
  P_\varepsilon u^{k,+}_\varepsilon&\to& u^{k}_0 \quad\text{in }\ \ H^1_0(\Omega)\text{-weak}\text{ as }E'\ni\varepsilon\to 0\label{eq3.15}\\
  P_\varepsilon u^{k,+}_\varepsilon&\to& u^{k}_0 \quad\text{in }\ \ L^2(\Omega)\text{ as }E'\ni\varepsilon\to 0\label{eq3.16}\\
\frac{\partial P_\varepsilon u^{k,+}_\varepsilon}{\partial
x_{j}}&\xrightarrow{2s} & \frac{\partial u_{0}^{k}}{\partial
x_{j}}+\frac{\partial u_{1}^{k}}{\partial y_{j}}\text{\ in
}L^{2}(\Omega)\text{ as }E'\ni\varepsilon\to 0\ (1\leq j\leq
N)\label{eq3.17}
\end{eqnarray}
where $(\lambda^{k}_0,u^{k}_0)\in \mathbb{R}\times
H^1_0(\Omega)$ is the $k^{th}$ eigencouple to the
 spectral problem
\begin{equation}\label{eq3.18}\left\{\begin{aligned}
 -\sum_{i,j=1}^N\frac{\partial}{\partial x_i}\left(\frac{1}{M_{S}(\rho)}q_{ij}\frac{\partial u_0}{\partial x_j}\right)
 &=\lambda_0 u_0\quad \text{in }\Omega\\
 u_0&=0\quad \text{ on  } \partial\Omega\\
 \int_{\Omega}|u_0|^2dx&=\frac{1}{M_{S}(\rho)},
\end{aligned}\right.
\end{equation}
and where $u_1^{k}\in L^2(\Omega;H^{1,*}_{per}(Y))$. Moreover, for almost every $x\in\Omega$ the following hold true:\\
\textbf{(i)} \ The restriction to $Y^*$ of $u_1^{k}(x)$ is the solution to the variational problem
\begin{equation}\label{eq3.21}\left\{\begin{aligned}
  &u_1^{k}(x)\in  H^1_{per}(Y^*)/\mathbb{R}\\&
  a(u_1^{k}(x),v)=-\sum_{i,j=1}^N \frac{\partial u^{k}_0}{\partial
  x_j}  \int_{Y^*}a_{ij}(y)\frac{\partial v}{\partial y_i}dy \\&
   \forall v\in  H^1_{per}(Y^*)/\mathbb{R};
   \end{aligned}\right.
\end{equation}
\textbf{(ii)} \ We have
\begin{equation}\label{eq3.22}
    u_1^{k}(x,y)=-\sum_{j=1}^N\frac{\partial u^{k}_0}{\partial
   x_j}(x)\chi^j(y)\qquad\text{a.e. in } (x,y)\in\Omega\times Y^*,
\end{equation}
where $\chi^j\ (1\leq j\leq N)$ is the solution to the cell problem (\ref{eq3.9}).
\end{theorem}
\begin{proof}We present only the outlines since this proof is similar but less technical to that
of the case $M_{S}(\rho)=0$.

Fix $k\geq 1$. By means of the minimax principle, as in
\cite{Vanni}, one easily proves the existence of a constant $C$
independent of $\varepsilon$ such that
$\frac{1}{\varepsilon}\lambda_\varepsilon^{k,+}<C$. Clearly, for fixed
$E\ni\varepsilon>0$, $u^{k,+}_\varepsilon$ lies in $V_\varepsilon$,
and
\begin{equation}\label{eq3.241}
\sum_{i,j=1}^N\int_{\Omega^\varepsilon}
a_{ij}(\frac{x}{\varepsilon})\frac{\partial
u^{k,+}_\varepsilon}{\partial x_j}
 \frac{\partial v}{\partial
 x_i}dx=\left(\frac{1}{\varepsilon}\lambda^{k,+}_\varepsilon\right)\varepsilon\int_{S^\varepsilon}\rho(\frac{x}{\varepsilon}) u^{k,+}_\varepsilon
 v\,d\sigma_\varepsilon(x)
\end{equation}
for any $v\in V_\varepsilon$. Bear in mind that
$\varepsilon\int_{S^\varepsilon}\rho(\frac{x}{\varepsilon})(u^{k,+}_\varepsilon)^2
dx=1$ and choose $v=u^{k,+}_\varepsilon$ in (\ref{eq3.241}). The
boundedness of the sequence
$(\frac{1}{\varepsilon}\lambda^{k,+}_\varepsilon)_{\varepsilon\in E}$
 and the ellipticity assumption (\ref{eq1.2}) imply at once by means of Proposition \ref{p3.1} that the
 sequence $(P_\varepsilon u^{k,+}_\varepsilon)_{\varepsilon\in E}$ is bounded in $H^1_0(\Omega)$.
  Theorem \ref{t2.2} and Proposition \ref{p2.1} apply simultaneously and
gives us $ \textbf{u}^{k}=(u_0^{k},u_1^{k})\in \mathbb{F}^1_0 $ such
that for some $\lambda_0^{k}\in\mathbb{R}$ and some subsequence
$E'\subset E$ we have (\ref{eq3.14})-(\ref{eq3.17}), where
(\ref{eq3.16}) is a direct consequence of (\ref{eq3.15}) by the
Rellich-Kondrachov theorem. For fixed $\varepsilon\in E'$, let
$\Phi_\varepsilon$ be as in Lemma \ref{l3.2}. Multiplying both sides
of  the first equality in (\ref{eq1.1}) by $\Phi_\varepsilon$ and
integrating over $\Omega^\varepsilon$ leads us to the variational
$\varepsilon$-problem
\begin{equation}\label{eq3.24}
\sum_{i,j=1}^N\int_{\Omega^\varepsilon}
a_{ij}(\frac{x}{\varepsilon})\frac{\partial P_\varepsilon
u^{k,+}_\varepsilon}{\partial x_j}
 \frac{\partial \Phi_\varepsilon}{\partial
 x_i}dx=(\frac{1}{\varepsilon}\lambda^{k,+}_\varepsilon)\varepsilon\int_{S^\varepsilon}(P_\varepsilon u^{k,+}_\varepsilon)\rho(\frac{x}{\varepsilon})
 \Phi_\varepsilon\,d\sigma_\varepsilon(x).
\end{equation}
Sending $\varepsilon\in E'$ to $0$, keeping
(\ref{eq3.14})-(\ref{eq3.17}) and Lemma \ref{l3.2} in mind, we
obtain
\begin{equation*}
\sum_{i,j=1}^N\iint_{\Omega\times Y^*}a_{ij}(y)
   \left(\frac{\partial u_0^k}{\partial x_j}+
    \frac{\partial u_1^k}{\partial y_j}\right)\left(\frac{\partial \psi_0}{\partial x_i}+
    \frac{\partial \psi_1}{\partial
    y_i}\right)dxdy=\lambda^{k}_0\iint_{\Omega\times S}u^{k}_0\psi_0(x)\rho(y)dxd\sigma(y).
\end{equation*}
Therefore, $(\lambda^{k}_0,\textbf{u}^{k})\in\mathbb{R}\times
\mathbb{F}^1_0$ solves the following \textit{global homogenized
spectral problem}:
\begin{equation}\label{eq3.26}\left\{\begin{aligned}
  &\text{Find }(\lambda,\textbf{u})\in\mathbb{C}\times
\mathbb{F}^1_0 \text{ such that }\\&
   \sum_{i,j=1}^N\iint_{\Omega\times Y^*}a_{ij}(y)
   \left(\frac{\partial u_0}{\partial x_j}+
    \frac{\partial u_1}{\partial y_j}\right)\left(\frac{\partial \psi_0}{\partial x_i}+
    \frac{\partial \psi_1}{\partial
    y_i}\right)dxdy=\lambda M_{S}(\rho)\int_\Omega u_0\psi_0
    \,dx\\&
   \text{for all } \Phi\in
\mathbb{F}^1_0 .
   \end{aligned}\right.
\end{equation}
which leads to the macroscopic and microscopic problems
(\ref{eq3.18})-(\ref{eq3.21}) without any major difficulty. As
regards the normalization condition in (\ref{eq3.18}), we fix
$k,l\geq 1$ and recall that the following holds for any
$\varphi\in\mathcal{D}(\Omega)$ (Proposition \ref{p2.1})
\begin{equation}\label{eq3.271}
\lim_{E'\ni\varepsilon\to
0}\varepsilon\int_{S^\varepsilon}(P_\varepsilon
u_{\varepsilon}^{k,+})\varphi(x)\rho(\frac{x}{\varepsilon})d\sigma_\varepsilon(x)=
\iint_{\Omega\times S}u_0^{k}(x)\varphi(x)\rho(y)\,dxd\sigma(y).
\end{equation}
But (\ref{eq3.271}) still holds for any $\varphi\in H^1_0(\Omega)$.
This being so, we write
\begin{eqnarray}
  && \varepsilon\int_{S^\varepsilon}(P_\varepsilon
u_{\varepsilon}^{k,+})(P_\varepsilon
u_{\varepsilon}^{l,+})\rho(\frac{x}{\varepsilon})d\sigma_\varepsilon(x)-M_{S}(\rho)\int_{\Omega}u_0^k
u_0^l\,dx\nonumber\\&
=&\varepsilon\int_{S^\varepsilon}(P_\varepsilon
u_{\varepsilon}^{k,+})(P_\varepsilon
u_{\varepsilon}^{l,+}-u_0^l)\rho(\frac{x}{\varepsilon})d\sigma_\varepsilon(x)
+\varepsilon\int_{S^\varepsilon}(P_\varepsilon
u_{\varepsilon}^{k,+})u_0^l\rho(\frac{x}{\varepsilon})d\sigma_\varepsilon(x)\label{eq3.272}\qquad\\&&-M_{S}(\rho)\int_{\Omega}u_0^k
u_0^l\,dx\nonumber
\end{eqnarray}
According to (\ref{eq3.271}) the sum of the last two terms on the
right hand side of (\ref{eq3.272}) goes to zero with $\varepsilon\in
E'$. As the remaining term on the right hand side of (\ref{eq3.272})
is concerned, we make use of the H\"{o}lder inequality to get
\begin{eqnarray*}
&& \left|\varepsilon\int_{S^\varepsilon}(P_\varepsilon
u_{\varepsilon}^{k,+})(P_\varepsilon
u_{\varepsilon}^{l,+}-u_0^l)\rho(\frac{x}{\varepsilon})d\sigma_\varepsilon(x)\right|\\&&\leq
 \|\rho\|_{\infty}\left(\varepsilon\int_{S^\varepsilon}|P_\varepsilon
u_{\varepsilon}^{k,+}|^2d\sigma_\varepsilon(x)\right)^\frac{1}{2}\left(\varepsilon\int_{S^\varepsilon}|P_\varepsilon
u_{\varepsilon}^{l,+}-u_0^l|^2d\sigma_\varepsilon(x)\right)^\frac{1}{2}.
\end{eqnarray*}
Next the trace inequality (see e.g., \cite{Pastukhova}) yields
\begin{eqnarray}
&&\varepsilon\int_{S^\varepsilon}|P_\varepsilon
u_{\varepsilon}^{k,+}|^2d\sigma_\varepsilon(x)\leq c
\left(\int_{\Omega^\varepsilon}|P_\varepsilon
u_{\varepsilon}^{k,+}|^2dx+\varepsilon^2\int_{\Omega^\varepsilon}|D
(P_\varepsilon u_{\varepsilon}^{k,+})|^2dx\right)\label{eq3.274}\\
&& \varepsilon\int_{S^\varepsilon}|P_\varepsilon
u_{\varepsilon}^{l,+}-u_0^l|^2d\sigma_\varepsilon(x)\leq
c\left(\int_{\Omega^\varepsilon}|P_\varepsilon
u_{\varepsilon}^{l,+}-u_0^l|^2dx+\varepsilon^2\int_{\Omega^\varepsilon}|D
(P_\varepsilon
u_{\varepsilon}^{l,+}-u_0^l)|^2dx\right)\qquad\qquad\label{eq3.273},
\end{eqnarray}
for some  positive constant $c$ independent of $\varepsilon$. But
the right hand side of (\ref{eq3.274}) is bounded from above whereas
that of (\ref{eq3.273}) converges to zero with $\varepsilon\in E'$.
This concludes the proof.
\end{proof}

\begin{remark}\label{r3.2}
\begin{itemize}
\item The eigenfunctions $\{u_0^{k}\}_{k=1}^\infty$ are in fact orthonormalized as follows
    $$
       \int_\Omega
       u_0^{k}u_0^{l}dx=\frac{\delta_{k,l}}{M_{S}(\rho)}\quad
       k,l=1,2,3,\cdots
    $$
    \item If $\lambda_0^{k}$ is simple (this is the case for
$\lambda_0^{1}$), then
by Theorem \ref{t3.1}, $\lambda_\varepsilon^{k,+}$ is also simple, for small $\varepsilon$,
and we can choose the eigenfunctions $u_\varepsilon^{k,+}$ such that the convergence results (\ref{eq3.15})-(\ref{eq3.17}) hold for
the whole sequence $E$.
\item Replacing $\rho$ with $-\rho$ in (\ref{eq1.1}), Theorem
 \ref{t3.1} also applies to the negative part of the spectrum in the
 case $M_{S}(\rho)<0$.
\end{itemize}
\end{remark}

\subsubsection{Negative part of the spectrum}
We now investigate the negative part of the spectrum  $(\lambda_\varepsilon^{k,-}, u_\varepsilon^{k,-})_{\varepsilon\in E}$.
 Before we can do this we need a few preliminaries and stronger regularity hypotheses on $T$, $\rho$ and the coefficients $(a_{ij})_{i,j=1}^N$.
 We assume in this subsection that $\partial T$ is $C^{2,\delta}$ and $\rho$ and the coefficients
  $(a_{ij})_{i,j=1}^N$ are $\delta$-H\"{o}lder continuous ($0<\delta<1$).

The following spectral problem is well posed
\begin{equation} \label{eq3.39}
\left\{\begin{aligned} &\text{Find }
(\lambda,\theta)\in\mathbb{C}\times
H^1_{per}(Y^*)\\
&-\sum_{i,j=1}^N\frac{\partial}{\partial
y_j}\left(a_{ij}(y)\frac{\partial
\theta}{\partial y_i}\right)=0 \text{ in }\ \ Y^*\\
&\sum_{i,j=1}^N a_{ij}(y)\frac{\partial \theta}{\partial y_i}\nu_j=\lambda\rho(y) \theta(y)
\text{ on } S
\end{aligned}\right.
\end{equation}
and possesses a spectrum with similar properties to that of (\ref{eq1.1}), two infinite (one positive and another negative) sequences. We recall that since we have $M_{S}(\rho)>0$, problem (\ref{eq3.39}) admits a unique nontrivial eigenvalue having an eigenfunction with definite sign, the first negative one (see e.g., \cite{Torne}). In the sequel we will only make use of $(\lambda_1^-, \theta_1^-)$, the first negative eigencouple to (\ref{eq3.39}). After proper sign choice we assume
that
\begin{equation}\label{eq3.40}
    \theta_1^-(y)>0 \ \  \text{ in }\  y\in Y^*.
\end{equation}
We also recall that $\theta_1^-$ is $\delta$-H\"{o}lder continuous(see e.g., \cite{Trudinger}), hence can be extended to a function living in  $\mathcal{C}_{per}(Y)$ still denoted by $\theta_1^-$. Notice that we have
\begin{equation}\label{eq3.391}
    \int_{S}\rho(y)(\theta_1^{-}(y))^2\,d\sigma(y)<0,
\end{equation}
as is easily seen from the variational equality (keep the
ellipticity hypothesis (\ref{eq1.2}) in mind)
$$
\sum_{i,j=1}^N\int_{Y^*}a_{ij}(y)\frac{\partial \theta_1^-}{\partial y_j}\frac{\partial \theta_1^-}{\partial y_i}\,dy=
\lambda_1^-\int_{S}\rho(y)(\theta_1^-(y))^2\,d\sigma(y).
$$
Bear in mind that problem (\ref{eq3.39}) induces by a scaling
argument the following equalities:
\begin{equation} \label{eq3.40}
\left\{\begin{aligned} &-\sum_{i,j=1}^N\frac{\partial}{\partial
x_j}\left(a_{ij}(\frac{x}{\varepsilon})\frac{\partial
\theta^\varepsilon}{\partial x_i}\right)=0\quad\text{ in } Q^\varepsilon\\
&\sum_{i,j=1}^N a_{ij}(\frac{x}{\varepsilon})\frac{\partial
\theta^\varepsilon}{\partial
x_i}\nu_j(\frac{x}{\varepsilon})=\frac{1}{\varepsilon}\lambda\rho(\frac{x}{\varepsilon})
\theta(\frac{x}{\varepsilon})\quad \text{ on } \partial
Q^\varepsilon,
\end{aligned}\right.
\end{equation}
where $\theta^\varepsilon(x)=\theta(\frac{x}{\varepsilon})$.
However, $\theta^\varepsilon$ is not zero on $\partial \Omega$. We
now introduce the following Steklov spectral problem (with an indefinite
density function)
\begin{equation} \label{eq3.41}\left\{\begin{aligned}
\text{Find } (\xi_\varepsilon,v_\varepsilon)\in\mathbb{C}\times
V_\varepsilon&\\
-\sum_{i,j=1}^N\frac{\partial}{\partial
x_j}\left(\widetilde{a}_{ij}(\frac{x}{\varepsilon})\frac{\partial
v_\varepsilon}{\partial x_i}\right)&=0\quad \text{ in } \Omega^\varepsilon\\
\sum_{i,j=1}^N
\widetilde{a}_{ij}(\frac{x}{\varepsilon})\frac{\partial
v_\varepsilon}{\partial
x_i}\nu_j(\frac{x}{\varepsilon})&=\xi_\varepsilon\widetilde{\rho}(\frac{x}{\varepsilon})
v_\varepsilon \text{ on } \partial
T^\varepsilon\\v_\varepsilon(x)&=0  \text{ on } \partial \Omega.
\end{aligned}\right.
\end{equation}
with new spectral parameters $(\xi_\varepsilon,
v_\varepsilon)\in\mathbb{C}\times V_\varepsilon$, where
$\widetilde{a}_{ij}(y)=(\theta_1^-)^2(y)a_{ij}(y)$ and
$\widetilde{\rho}(y)=(\theta_1^-)^2(y)\rho(y)$. Notice that $\widetilde{a}_{ij}(y)\in L^\infty_{per}(Y)$ and
$\widetilde{\rho}(y)\in \mathcal{C}_{per}(Y)$. As
$0<c_-\leq \theta_1^-(y)\leq c^+<+\infty$ ($c_-, c^+\in\mathbb{R}$), the operator on the left hand side of
(\ref{eq3.41}) is uniformly elliptic and Theorem \ref{t3.1} applies to the negative part of the spectrum of
(\ref{eq3.41}) (see (\ref{eq3.391}) and Remark \ref{r3.2}). The effective spectral problem for (\ref{eq3.41}) reads
\begin{equation}\label{eq3.411}\left\{\begin{aligned}
 -\sum_{i,j=1}^N\frac{\partial}{\partial x_j}\left(\widetilde{q}_{ij}
 \frac{\partial v_0}{\partial x_i}\right)&=\xi_0 M_{S}(\widetilde{\rho}) v_0\quad \text{in
 }\Omega\\v_0&=0\quad \text{ on  } \partial\Omega\\\int_{\Omega}|v_0|^2dx&=\frac{-1}{M_{S}(\widetilde{\rho})}.
\end{aligned}\right.
\end{equation}
The effective coefficients $\{\widetilde{q}_{ij}\}_{1\leq i,j\leq N}$
 being defined as expected, i.e.,
\begin{equation}\label{eq3.42}
\widetilde{q}_{ij}=\int_{Y^*}\widetilde{a}_{ij}(y)dy-\sum_{l=1}^N\int_{Y^*}\widetilde{a}_{il}(y)\frac{\partial
\widetilde{\chi}^j_1}{\partial y_l}(y)dy,
\end{equation}
with $\widetilde{\chi}^l\in H^1_{per}(Y^*)/\mathbb{R}$  $(l=1,...,N)$ being the
solution to the following local problem
\begin{equation} \label{eq3.43}\left\{\begin{aligned}
&\widetilde{\chi}^l\in H^1_{per}(Y^*)/\mathbb{R}\\&
\sum_{i,j=1}^N\int_{Y^*}\widetilde{a}_{ij}(y)\frac{\partial
\widetilde{\chi}^l}{\partial y_j}\frac{\partial
v}{\partial
y_i}dy=\sum_{i=1}^N\int_{Y^*}\widetilde{a}_{il}(y)\frac{\partial
v}{\partial y_i}dy\\& \text{for all } v\in H^1_{per}(Y^*)/\mathbb{R}
.\end{aligned}\right.
\end{equation}
We will use the following notation in the sequel:
\begin{equation*}
\widetilde{a}(u,v)=\sum_{i,j=1}^N\int_{Y^*}\widetilde{a}_{ij}(y)\frac{\partial
u}{\partial y_j}\frac{\partial
v}{\partial
y_i}dy\qquad \left(u,v\in H^1_{per}(Y^*)/\mathbb{R}\right).
\end{equation*}
Notice that the spectrum of (\ref{eq3.411}) is as follows
$$
0>\xi_0^1>\xi_0^2\geq \xi_0^3\geq \cdots \geq\xi_0^j \geq \cdots   \to -\infty \text { as } j\to\infty.
$$
Making use of (\ref{eq3.40}) when following the same line of reasoning as in \cite[Lemma 6.1]{Vanni}, we obtain
that the negative spectral parameters of problems (\ref{eq1.1}) and (\ref{eq3.41}) verify:
\begin{equation}\label{eq3.414}
u_\varepsilon^{k,-}=(\theta_1^-)^\varepsilon v_\varepsilon^{k,-}\quad (\varepsilon\in E,\ k=1,2\cdots)
\end{equation}
and
\begin{equation}\label{eq3.415}
\lambda_\varepsilon^{k,-}=\frac{1}{\varepsilon}\lambda^-_1+\xi_\varepsilon^{k,-} + o(1) \quad (\varepsilon\in E,\ k=1,2\cdots).
\end{equation}
The presence of the term $o(1)$ is due to integrals over
$\Omega^\varepsilon\setminus Q^\varepsilon$, which converge to zero with $\varepsilon$, remember that (\ref{eq3.40})
 holds in $Q^\varepsilon$ but not
$\Omega^\varepsilon$. This trick, known as "factorization principle" was introduced by Vaninathan\cite{Vanni} and has
been used in many other works on averaging, see e.g., \cite{AP, Kozlov, Nazarov2} just to cite a few. As will be seen
below, the sequence $ (\xi_\varepsilon^{k,-})_{\varepsilon\in E} $ is bounded in $\mathbb{R}$. In another words,
$\lambda_\varepsilon^{k,- }$ is of order $1/\varepsilon$ and tends to $-\infty$ as $\varepsilon$ goes to zero. It is now clear why the limiting behavior of negative eigencouples is not straightforward as that of positive ones and requires further investigations, which have just been made.

Indeed, as the reader might be guessing now, the suitable orthonormalization condition for (\ref{eq3.41}) is
\begin{equation}
    \varepsilon\int_{S^\varepsilon}\widetilde{\rho}(\frac{x}{\varepsilon})v_\varepsilon^{k,-}v_\varepsilon^{l,-}\,d\sigma_\varepsilon(x)=-\delta_{k,l}\quad
    k,l=1,2,\cdots
\end{equation}
which by means of $(\ref{eq3.414})$ is equivalent to
\begin{equation}\label{eq3.412}
    \varepsilon\int_{S^\varepsilon}\rho(\frac{x}{\varepsilon})u_\varepsilon^{k,-}u_\varepsilon^{l,-}\,d\sigma_\varepsilon(x)=-\delta_{k,l}\quad
    k,l=1,2,\cdots
\end{equation}
We may now state the homogenization theorem for the negative part of the spectrum of (\ref{eq1.1}).

\begin{theorem}\label{t3.2}
For each $k\geq 1$ and each $\varepsilon\in E$, let
$(\lambda^{k,-}_\varepsilon,u^{k,-}_\varepsilon)$ be the $k^{th}$
negative eigencouple to (\ref{eq1.1}) with $M_{S}(\rho)>0$ and
(\ref{eq3.412}). Then, there exists a subsequence $E'$ of $E$ such
that
\begin{eqnarray}
 \frac{\lambda^{k,-}_\varepsilon}{\varepsilon}-\frac{\lambda_1^-}{\varepsilon^2} &\to&  \xi^k_0\quad\text{in }\ \mathbb{R}\ \text{ as }E\ni\varepsilon\to 0\label{eq3.44}\\
  P_\varepsilon v^{k,-}_\varepsilon&\to& v^{k}_0 \quad\text{in }\ \ H^1_0(\Omega)\text{-weak}\text{ as }E'\ni\varepsilon\to 0\label{eq3.45}\\
  P_\varepsilon  v^{k,-}_\varepsilon&\to& v^{k}_0 \quad\text{in }\ \ L^2(\Omega)\text{ as }E'\ni\varepsilon\to 0\label{eq3.46}\\
\frac{\partial P_\varepsilon  v^{k,-}_\varepsilon}{\partial
x_{j}}&\xrightarrow{2s} & \frac{\partial v_{0}^{k}}{\partial
x_{j}}+\frac{\partial v_{1}^{k}}{\partial y_{j}}\text{\ in
}L^{2}(\Omega)\text{ as }E'\ni\varepsilon\to 0\ (1\leq j\leq
N)\label{eq3.47}
\end{eqnarray}
where $(\xi^{k}_0,v^{k}_0)\in \mathbb{R}\times H^1_0(\Omega)$ is the
 $k^{th}$ eigencouple to the
 spectral problem
\begin{equation}\label{eq3.48}\left\{\begin{aligned}
 -\sum_{i,j=1}^N\frac{\partial}{\partial x_i}\left(\frac{1}{M_{S}(\widetilde{\rho})}\widetilde{q}_{ij}
 \frac{\partial v_0}{\partial x_j}\right)&=\xi_0 v_0\quad \text{in
 }\Omega\\v_0&=0\quad \text{ on  } \partial\Omega\\\int_{\Omega}|v_0|^2\,dx&=\frac{-1}{M_{S}(\widetilde{\rho})},
\end{aligned}\right.
\end{equation}
and where $v_1^{k}\in L^2(\Omega;H^{1,*}_{per}(Y))$. Moreover, for almost every $x\in\Omega$ the following hold true:\\
\textbf{(i)} \ The restriction to $Y^*$ of $v_1^{k}(x)$ is the solution to the variational problem
\begin{equation}\label{eq3.49}\left\{\begin{aligned}
  &v_1^{k}(x)\in  H^1_{per}(Y^*)/\mathbb{R}\\&
  \widetilde{a}(v_1^{k}(x),u)=-\sum_{i,j=1}^N \frac{\partial v^{k}_0}{\partial
  x_j}  \int_{Y^*}\widetilde{a}_{ij}(y)\frac{\partial u}{\partial y_i}dy \\&
   \forall u\in  H^1_{per}(Y^*)/\mathbb{R};
   \end{aligned}\right.
\end{equation}
\textbf{(ii)} \ We have
\begin{equation}\label{eq3.50}
    v_1^{k}(x,y)=-\sum_{j=1}^N\frac{\partial v^{k}_0}{\partial
   x_j}(x)\widetilde{\chi}^j(y)\qquad\text{a.e. in } (x,y)\in\Omega\times Y^*,
\end{equation}
where $\widetilde{\chi}^j\ (1\leq j\leq N)$ is the solution  to the cell problem (\ref{eq3.43}).
\end{theorem}
\begin{remark}
\begin{itemize}
\item The eigenfunctions $\{v_0^{k}\}_{k=1}^\infty$ are orthonormalized by
    $$
       \int_\Omega
       v_0^{k}v_0^{l}dx=\frac{-\delta_{k,l}}{M_{S}(\widetilde{\rho})}\quad
       k,l=1,2,3,\cdots
    $$
    \item If $\xi_0^{k}$ is simple (this is the case for
$\xi_0^{1}$), then
by Theorem \ref{t3.2}, $\lambda_\varepsilon^{k,-}$ is also simple for small $\varepsilon$,
and we can choose the `eigenfunctions' $v_\varepsilon^{k,-}$ such that the convergence results (\ref{eq3.45})-(\ref{eq3.47}) hold for
the whole sequence $E$.
\item Replacing $\rho$ with $-\rho$ in (\ref{eq1.1}), Theorem
 \ref{t3.2} adapts to the positive part of the spectrum in the
 case $M_{S}(\rho)<0$.
\end{itemize}
\end{remark}

\subsection{The case $M_{S}(\rho)=0$}
We prove an homogenization result for both the positive part and the
negative part of the spectrum simultaneously. We assume in this case
that the eigenfunctions are orthonormalized as follows
\begin{equation}\label{eq3.22}
   \int_{S^\varepsilon}\rho(\frac{x}{\varepsilon})u_\varepsilon^{k,\pm}u_\varepsilon^{l,\pm}\,d\sigma_\varepsilon(x)=\pm\delta_{k,l}\quad
    k,l=1,2,\cdots
\end{equation}
Let $\chi^0$ be the solution to $(\ref{eq3.91})$ and put
\begin{equation}\label{eq3.23}
    \nu^2=\sum_{i,j=1}^N\int_{Y^*}a_{ij}(y)\frac{\partial \chi^0}{\partial y_j}\frac{\partial \chi^0}{\partial
    y_i}dy.
\end{equation}
Indeed, the right hand side of (\ref{eq3.23}) is positive. We recall
that the following spectral problem for a quadratic operator pencil
with respect to $\nu$,
\begin{equation}\label{eq3.24}
\left\{\begin{aligned} -\sum_{i,j=1}^N\frac{\partial}{\partial
x_j}\left(q_{ij}\frac{\partial
u_0}{\partial x_i}\right)&=\lambda_0^2\nu^2 u_0\text{ in } \Omega\\
u_0&=0 \text{ on } \partial \Omega,
\end{aligned}\right.
\end{equation}
has a spectrum consisting of two infinite sequences
$$
0<\lambda_0^{1,+}< \lambda_0^{2,+}\leq \cdots \leq
\lambda_0^{k,+}\leq \dots,\quad \lim_{n\to
+\infty}\lambda_0^{k,+}=+\infty
$$
and
$$
0>\lambda_0^{1,-}> \lambda_0^{2,-}\geq \cdots \geq
\lambda_0^{k,-}\geq \dots,\quad \lim_{n\to
+\infty}\lambda_0^{k,-}=-\infty.
$$
with $\lambda^{k,+}_0=-\lambda^{k,-}_0\ \ k=1,2,\cdots$ and with the
corresponding eigenfunctions $u_0^{k,+}=u_0^{k,-}$. We note by
passing that $\lambda^{1,+}_0$ and $\lambda^{1,-}_0$ are simple. We
are now in a position to state the homogenization result in the
present case.

\begin{theorem}\label{t3.3}
For each $k\geq 1$ and each $\varepsilon\in E$, let
$(\lambda^{k,\pm}_\varepsilon,u^{k,\pm}_\varepsilon)$ be the
$(k,\pm)^{th}$ eigencouple to (\ref{eq1.1}) with $M_{S}(\rho)=0$
and (\ref{eq3.22}). Then, there exists a subsequence $E'$ of $E$
such that
\begin{eqnarray}
 \lambda^{k,\pm}_\varepsilon &\to&  \lambda^{k,\pm}_0\quad\text{in }\ \mathbb{R}\ \text{ as }E\ni\varepsilon\to 0\label{eq3.25}\\
  P_\varepsilon u^{k,\pm}_\varepsilon&\to& u^{k,\pm}_0 \quad\text{in }\ \ H^1_0(\Omega)\text{-weak}\text{ as }E'\ni\varepsilon\to 0\label{eq3.26}\\
  P_\varepsilon u^{k,\pm}_\varepsilon&\to& u^{k,\pm}_0 \quad\text{in }\ \ L^2(\Omega)\text{ as }E'\ni\varepsilon\to 0\label{eq3.27}\\
\frac{\partial P_\varepsilon u^{k,\pm}_\varepsilon}{\partial
x_{j}}&\xrightarrow{2s} & \frac{\partial u_{0}^{k,\pm}}{\partial
x_{j}}+\frac{\partial u_{1}^{k,\pm}}{\partial y_{j}}\text{\ in
}L^{2}(\Omega)\text{ as }E'\ni\varepsilon\to 0\ (1\leq j\leq
N)\label{eq3.28}
\end{eqnarray}
where $(\lambda^{k,\pm}_0,u^{k,\pm}_0)\in \mathbb{R}\times
H^1_0(\Omega)$ is the $(k,\pm)^{th}$ eigencouple to the following
 spectral problem for a quadratic operator pencil with respect to $\nu$,
\begin{equation}\label{eq3.29}
\left\{\begin{aligned} -\sum_{i,j=1}^N\frac{\partial}{\partial
x_i}\left(q_{ij}\frac{\partial
u_0}{\partial x_j}\right)&=\lambda_0^2\nu^2 u_0\text{ in } \Omega\\
u_0&=0 \text{ on } \partial \Omega,
\end{aligned}\right.
\end{equation}
and where $u_1^{k,\pm}\in L^2(\Omega;H^{1,*}_{per}(Y))$.  We have the following normalization condition
\begin{equation}\label{eq2.291}
    \int_{\Omega}|u_0^{k,\pm}|^2\,dx=\frac{\pm 1}{2\lambda_0^{k,\pm} \nu^2}\qquad k=1,2,\cdots
\end{equation}
Moreover, for almost every $x\in\Omega$ the following hold true:\\
\textbf{(i)} \ The restriction to $Y^*$ of $u_1^{k,\pm}(x)$ is the solution to the variational problem
\begin{equation}\label{eq3.30}\left\{\begin{aligned}
  &u_1^{k,\pm}(x)\in  H^1_{per}(Y^*)/\mathbb{R}\\&
  a(u_1^{k,\pm}(x),v)=\lambda^{k,\pm}_0 u_0^{k,\pm}(x)\int_{S}\rho(y)v(y)\,d\sigma(y)-\sum_{i,j=1}^N \frac{\partial u^{k,\pm}_0}{\partial
  x_j}(x)  \int_{Y^*}a_{ij}(y)\frac{\partial v}{\partial y_i}\,dy \\&
   \forall v\in  H^1_{per}(Y^*)/\mathbb{R};
   \end{aligned}\right.
\end{equation}
\textbf{(ii)} \ We have
\begin{equation}\label{eq3.31}
    u_1^{k,\pm}(x,y)=\lambda^{k,\pm}_0 u_0^{k,\pm}(x)\chi^0(y)-\sum_{j=1}^N\frac{\partial u^{k,\pm}_0}{\partial
   x_j}(x)\chi^j(y)\qquad\text{a.e. in } (x,y)\in\Omega\times Y^*,
\end{equation}
where $\chi^j\ (1\leq j\leq N)$ and $\chi^0$ are the solutions to the cell problems (\ref{eq3.9}) and
(\ref{eq3.91}), respectively.
\end{theorem}

\begin{proof}
Fix $k\geq 1$, using the minimax principle, as in \cite{Vanni}, we
get a constant $C$ independent of $\varepsilon$ such that
$|\lambda_\varepsilon^{k,\pm}|<C$. We have
$u^{k,\pm}_\varepsilon\in V_\varepsilon$ and
\begin{equation}\label{eq3.32}
\sum_{i,j=1}^N\int_{\Omega^\varepsilon}
a_{ij}(\frac{x}{\varepsilon})\frac{\partial
  u^{k,\pm}_\varepsilon}{\partial x_j}
 \frac{\partial v}{\partial
 x_i}dx=\lambda^{k,\pm}_\varepsilon\int_{S^\varepsilon}  \rho(\frac{x}{\varepsilon})u^{k,\pm}_\varepsilon
 v\,d\sigma_\varepsilon(x)
\end{equation}
for any $v\in V_\varepsilon$. Bear in mind that
$\int_{S^\varepsilon}\rho(\frac{x}{\varepsilon})(u^{k,\pm}_\varepsilon)^2
\,d\sigma_\varepsilon(x)=\pm 1$ and choose $v=u^{k,\pm}_\varepsilon$ in (\ref{eq3.32}).
The boundedness of the sequence $(\lambda^{k,\pm}_\varepsilon)_{\varepsilon\in E}$
 and the ellipticity assumption (\ref{eq1.2}) imply at once by means of Proposition \ref{p3.1} that the
 sequence $(P_\varepsilon u^{k,\pm}_\varepsilon)_{\varepsilon\in E}$ is bounded in $H^1_0(\Omega)$.
  Theorem \ref{t2.2} and Proposition \ref{p2.1} apply simultaneously and
gives us $ \textbf{u}^{k,\pm}=(u_0^{k,\pm},u_1^{k,\pm})\in
\mathbb{F}^1_0 $ such that for some $\lambda_0^{k,\pm}\in\mathbb{R}$
and some subsequence $E'\subset E$ we have
(\ref{eq3.25})-(\ref{eq3.28}), where (\ref{eq3.27}) is a direct
consequence of (\ref{eq3.26}) by the Rellich-Kondrachov theorem. For
fixed $\varepsilon\in E'$, let $\Phi_\varepsilon$ be as in Lemma
\ref{l3.2}. Multiplying both sides of  the first equality in
(\ref{eq1.1}) by $\Phi_\varepsilon$ and integrating over
$\Omega^\varepsilon$ leads us to the variational
$\varepsilon$-problem
\begin{equation*}
\sum_{i,j=1}^N\int_{\Omega^\varepsilon}
a_{ij}(\frac{x}{\varepsilon})\frac{\partial P_\varepsilon
u^{k,\pm}_\varepsilon}{\partial x_j}
 \frac{\partial \Phi_\varepsilon}{\partial
 x_i}\,dx=\lambda^{k,\pm}_\varepsilon\int_{S^\varepsilon}(P_\varepsilon u^{k,\pm}_\varepsilon)\rho(\frac{x}{\varepsilon})
\Phi_\varepsilon\,d\sigma_\varepsilon(x).
\end{equation*}
Sending $\varepsilon\in E'$ to $0$, keeping
(\ref{eq3.25})-(\ref{eq3.28}) and Lemma \ref{l3.2} in mind, we
obtain
\begin{equation}\label{eq3.33}
a_\Omega(\mathbf{u}^{k,\pm},\Phi)=\lambda^{k,\pm}_0\iint_{\Omega\times
S}\left(u_1^{k,\pm}(x,y)\psi_0(x)\rho(y)+
    u_0^{k,\pm}\psi_1(x,y)\rho(y)\right)\,dxd\sigma(y)
    \end{equation}
The right-hand side follows as explained below. we have

\begin{eqnarray*}
  \int_{S^\varepsilon}(P_\varepsilon
u^{k,\pm}_\varepsilon)\rho(\frac{x}{\varepsilon})
\Phi_\varepsilon\,d\sigma_\varepsilon(x)&=&\int_{S^\varepsilon}(P_\varepsilon
u^{k,\pm}_\varepsilon)\psi_0(x)
 \rho(\frac{x}{\varepsilon})
 \,d\sigma_\varepsilon(x) \\
   &+&\varepsilon\int_{S^\varepsilon}(P_\varepsilon u^{k,\pm}_\varepsilon)\psi_1(x,\frac{x}{\varepsilon})
 \rho(\frac{x}{\varepsilon}) \,d\sigma_\varepsilon(x).
\end{eqnarray*}
On the one hand we have
$$
\lim_{E'\ni\varepsilon\to 0}\varepsilon\int_{S^\varepsilon}(P_\varepsilon
u^{k,\pm}_\varepsilon)\psi_1(x,\frac{x}{\varepsilon})
 \rho(\frac{x}{\varepsilon})
 \,dx=\iint_{\Omega\times
S} u_0^{k,\pm}\psi_1(x,y)\rho(y)\,dxd\sigma(y).
$$
On the other hand, owing to Lemma \ref{l2.1}, the following holds:
$$
\lim_{E'\ni\varepsilon\to 0}\int_{S^\varepsilon}
(P_\varepsilon u^{k,\pm}_\varepsilon)\psi_0(x)
 \rho(\frac{x}{\varepsilon})
 \,d\sigma_\varepsilon(x)=\iint_{\Omega\times
S}u_1^{k,\pm}(x,y)\psi_0(x)\rho(y)\,dxd\sigma(y).
$$
We have just proved that
$(\lambda^{k,\pm}_0,\textbf{u}^{k,\pm})\in\mathbb{R}\times
\mathbb{F}^1_0$ solves the following \textit{global homogenized
spectral problem}:
\begin{equation}\label{eq3.34}\left\{\begin{aligned}
  &\text{Find }(\lambda,\textbf{u})\in\mathbb{C}\times
\mathbb{F}^1_0 \text{ such that }\\&
  a_\Omega(\mathbf{u},\Phi)=\lambda\iint_{\Omega\times
S}\left(u_1(x,y)\psi_0(x)+
    u_0(x)\psi_1(x,y)\right)\rho(y)\,dxd\sigma(y) \\&
   \text{for all } \Phi\in
\mathbb{F}^1_0 .
   \end{aligned}\right.
\end{equation}
To prove (i), choose $\Phi=(\psi_0,\psi_1)$ in (\ref{eq3.33})  such
that $\psi_{0}=0$ and $\psi_{1}=\varphi\otimes v_1$, where
$\varphi\in\mathcal{D}(\Omega)$ and $v_1\in H^1_{per}(Y^*)/\mathbb{R}$ to get
\[
\int_{\Omega}\varphi(x)\left[\sum_{i,j=1}^N\int_{Y^*}a_{ij}(y)\left(\frac{\partial
u^{k,\pm}_0}{\partial x_j}+\frac{\partial u^{k,\pm}_1}{\partial y_j}
\right)\frac{\partial v_1}{\partial
y_i}dy\right]dx=\int_\Omega\varphi(x)\left[\lambda_0^{k,\pm}u_0^{k,\pm}(x)\int_{S}v_1(y)\rho(y)\,d\sigma(y)\right]dx
\]
Hence by the arbitrariness of $\varphi$, we have a.e. in $\Omega$
\[
\sum_{i,j=1}^N\int_{Y^*}a_{ij}(y)\left(\frac{\partial
u^{k,\pm}_0}{\partial x_j}+\frac{\partial u^{k,\pm}_1}{\partial y_j}
\right)\frac{\partial v_1}{\partial
y_i}dy=\lambda_0^{k,\pm}u_0^{k,\pm}(x)\int_{S}v_1(y)\rho(y)d\sigma(y)
\]
for any  $v_1$ in  $H^1_{per}(Y^*)/\mathbb{R}$, which is nothing but (\ref{eq3.30}).

Fix $x\in\overline{\Omega}$, multiply both sides of (\ref{eq3.9}) by
$-\frac{\partial u_0^{k,\pm}}{\partial x_j}(x)$ and sum over $1\leq
j\leq N$. Adding side by side to the resulting equality that
obtained after multiplying both sides of (\ref{eq3.91}) by
$\lambda_0^{k,\pm}u_0^{k,\pm}(x)$, we realize that
$z(x)=-\sum_{j=1}^N\frac{\partial u_0^{k,\pm}}{\partial
x_j}(x)\chi^j(y)+\lambda_0^{k,\pm}u_0^{k,\pm}(x)\chi^0(y)$ solves
(\ref{eq3.30}). Hence
\begin{equation}\label{eq3.36}
     u_1^{k,\pm}(x,y)=\lambda^{k,\pm}_0 u_0^{k,\pm}(x)\chi^0(y)-\sum_{j=1}^N\frac{\partial u^{k,\pm}_0}{\partial
   x_j}(x)\chi^j(y)\ \ \ a.e.\ \text{in }\ \Omega\times Y^*,
\end{equation}
by uniqueness of the solution to (\ref{eq3.30}). Thus (\ref{eq3.31}).
But (\ref{eq3.36}) still holds almost everywhere in $(x,y)\in\Omega\times S$ as $S$ is of class $\mathcal{C}^1$.
Considering now
$\Phi=(\psi_0,\psi_1)$ in (\ref{eq3.33})  such that
$\psi_{0}\in\mathcal{D}(\Omega)$ and $\psi_{1}=0$ we get
$$
\sum_{i,j=1}^N\iint_{\Omega\times Y^*}a_{ij}(y)\left(\frac{\partial
u_0^{k,\pm}}{\partial x_j}+\frac{\partial u_1^{k,\pm}}{\partial
y_j}\right)\frac{\partial \psi_0}{\partial
x_i}\,dxdy=\lambda_0^{k,\pm}\iint_{\Omega\times
S}u_1^{k,\pm}(x,y)\rho(y)\psi_0(x)\,dxd\sigma(y),
$$
which by means of (\ref{eq3.36}) leads to
\begin{eqnarray}
    &&\sum_{i,j=1}^N  \int_{\Omega}q_{ij}\frac{\partial
u_0^{k,\pm}}{\partial x_j}\frac{\partial \psi_0}{\partial
x_i}dx+\lambda_0^{k,\pm}\sum_{i,j=1}^N\int_\Omega
u_0^{k,\pm}(x)\frac{\partial \psi_0}{\partial
x_i}\left(\int_{Y^*}a_{ij}(y)\frac{\partial \chi^0}{\partial y_j
}(y)dy\right)dx\nonumber\\
&&=-\lambda_0^{k,\pm}\sum_{j=1}^N\int_\Omega \frac{\partial
u_0^{k,\pm}}{\partial x_j}
\psi_0(x)\left(\int_{S}\rho(y) \chi^j(y)\,d\sigma(y)\right)dx\label{eq3.37}\\
&&+(\lambda_0^{k,\pm})^2\int_\Omega
u_0^{k,\pm}(x)\psi_0(x)\left(\int_{S}\rho(y)\chi^0(y)\,d\sigma(y)\right)dx\nonumber.
\end{eqnarray}
Choosing $\chi^l\ (1\leq l\leq N)$ as test function in
(\ref{eq3.91}) and $\chi^0$ as test function in (\ref{eq3.9}) we
observe that
$$
\sum_{j=1}^N\int_{Y^*}a_{lj}(y)\frac{\partial \chi^0}{\partial y_j
}(y)dy=\int_{S}\rho(y) \chi^l(y)\,d\sigma(y)=a(\chi^l,\chi^0)\quad
(l=1,\cdots N).
$$
Thus, in (\ref{eq3.37}), the second term in the left hand side is
equal to the first one in the right hand side. This leaves us with
\begin{equation}\label{eq3.38}
\int_{\Omega}q_{ij}\frac{\partial u_0^{k,\pm}}{\partial
x_j}\frac{\partial \psi_0}{\partial
x_i}dx=(\lambda_0^{k,\pm})^2\int_\Omega
u_0^{k,\pm}(x)\psi_0(x)dx\left(\int_{S}\rho(y)\chi^0(y)\,d\sigma(y)\right).
\end{equation}
Choosing $\chi^0$ as test function in (\ref{eq3.91}) reveals that
$$
\int_{S}\rho(y)\chi^0(y)\,d\sigma(y)=a(\chi^0,\chi^0)=\nu^2.
$$
Hence
$$
\sum_{i,j=1}^N\int_{\Omega}q_{ij}\frac{\partial
u_0^{k,\pm}}{\partial x_j}\frac{\partial \psi_0}{\partial
x_i}dx=(\lambda_0^{k,\pm})^2\nu^2\int_\Omega
u_0^{k,\pm}(x)\psi_0(x)dx,
$$
and
$$
-\sum_{i,j=1}^N\frac{\partial}{\partial
x_i}\left(q_{ij}\frac{\partial u_0^{k,\pm}}{\partial
x_j}(x)\right)=(\lambda_0^{k,\pm})^2\nu^2 u_0^{k,\pm}(x)\text{ in }
\Omega.
$$
Thus, the convergence (\ref{eq3.25}) holds for the whole sequence
$E$. We now address (\ref{eq2.291}). Fix  $k,l\geq 1$ and let
$\vartheta\in H^1_{per}(Y^*)/\mathbb{R}$ be the solution to
(\ref{eq3.998}) where $\theta$ is replaced with  our density
function $\rho$. As in (\ref{eq3.999}), we transform the surface
integral into a volume integral
\begin{eqnarray}
  \int_{S^\varepsilon}(P_\varepsilon u_{\varepsilon}^{k,\pm})(P_\varepsilon u_{\varepsilon}^{l,\pm})\rho
  (\frac{x}{\varepsilon})d\sigma_\varepsilon(x) &=& \int_{\Omega^\varepsilon}(P_\varepsilon u_{\varepsilon}^{k,\pm})
  D_x(P_\varepsilon u_{\varepsilon}^{l,\pm})\cdot D_y\vartheta(\frac{x}{\varepsilon}) dx \nonumber\\
   &+& \int_{\Omega^\varepsilon}D_x(P_\varepsilon u_{\varepsilon}^{k,\pm})(P_\varepsilon u_{\varepsilon}^{l,\pm})\cdot
   D_y\vartheta(\frac{x}{\varepsilon}) dx.\qquad\label{eq3.266}
\end{eqnarray}
A limit passage in (\ref{eq3.266}) as $E'\ni\varepsilon\to 0$ yields
\begin{eqnarray*}
  && \lim_{E'\ni\varepsilon\to 0} \int_{S^\varepsilon}(P_\varepsilon u_{\varepsilon}^{k,\pm})(P_\varepsilon u_{\varepsilon}^{l,\pm})\rho(\frac{x}{\varepsilon})d\sigma_\varepsilon(x) \qquad\qquad\qquad\qquad \\
   &=&\iint_{\Omega\times Y^*}u_0^{k,\pm}(D_x u_0^{l,\pm}+D_y u_1^{l,\pm})\cdot D_y\vartheta dxdy+ \iint_{\Omega\times Y^*}(D_x u_0^{k,\pm}+D_y u_1^{k,\pm})u_0^{l,\pm}\cdot D_y\vartheta dxdy\qquad\qquad\qquad\qquad \\
   &=& \int_{\Omega}u_0^{k,\pm}\left(\int_{Y^*}D_y u_1^{l,\pm}(x,y)\cdot D_y\vartheta(y)dy\right)dx +\int_{\Omega}u_0^{l,\pm}\left(\int_{Y^*}D_y u_1^{k,\pm}(x,y)\cdot D_y\vartheta(y)dy\right)dx\\
   &=&\iint_{\Omega\times S}u_0^{k,\pm}(x) u_1^{l,\pm}(x,y)\rho(y)dxd\sigma(y)+\iint_{\Omega\times S}u_0^{l,\pm}(x) u_1^{k,\pm}(x,y)\rho(y)dxd\sigma(y)  \\
   &=& \lambda_0^{l,\pm} \nu^2\int_{\Omega}u_0^{k,\pm}(x)u_0^{l,\pm}(x)dx+\lambda_0^{k,\pm} \nu^2\int_{\Omega}u_0^{l,\pm}(x)u_0^{k,\pm}(x)dx\\
   &=&( \lambda_0^{k,\pm}+\lambda_0^{l,\pm})\nu^2\int_{\Omega}u_0^{k,\pm}(x)u_0^{l,\pm}(x)dx.
\end{eqnarray*}
Where after the limit passage, we used the integration by part
formula, then the weak formulation of (\ref{eq3.998})
 and finally (\ref{eq3.31}) and integration by part. If $k=l$, the above limit passage and (\ref{eq3.22}) lead to the
 desired result, (\ref{eq2.291}), completing thereby the proof.
\end{proof}
\begin{remark}
\begin{itemize}
      \item The eigenfunctions $\{u_0^{k,\pm}\}_{k=1}^\infty$ are in fact orthonormalized as follows
      $$
       \int_{\Omega}u_0^{l,\pm}(x)u_0^{k,\pm}(x)dx=\frac{\pm\delta_{k,l}}{\nu^2(\lambda_0^{l,\pm}+\lambda_0^{k,\pm})}\qquad k,l =1,2,\cdots\nonumber
      $$
     \item If $\lambda_0^{k,\pm}$ is simple (this is the case for
$\lambda_0^{1,\pm}$), then by Theorem \ref{t3.3},
$\lambda_\varepsilon^{k,\pm}$ is also simple, for small
$\varepsilon$, and we can choose the eigenfunctions
$u_\varepsilon^{k,\pm}$ such that the convergence results
(\ref{eq3.26})-(\ref{eq3.28}) hold for the whole sequence $E$.
\end{itemize}
\end{remark}

\subsection*{Final Remark}
After this paper was completed (see \cite{douanla4}) and submitted, we learned about an independent and similar work \cite{PN}.

\subsection*{Acknowledgments}
The author is grateful to Dr. Jean Louis Woukeng for helpful
discussions.

\end{document}